\documentclass[dvips, preprint]{imsart}

\RequirePackage[OT1]{fontenc}
\RequirePackage{amsthm,amsmath,natbib}

\usepackage{amssymb, epsfig, amsfonts}

\startlocaldefs
\numberwithin{equation}{section}
\theoremstyle{plain}

\newtheorem{prop}{Proposition}[section]
\newtheorem{lem}{Lemma}[section]
\newtheorem{cor}{Corollary}[section]
\endlocaldefs

\begin{document}

%\baselineskip = 1.5\baselineskip

\begin{frontmatter}

\title{A general semiparametric Z-estimation approach for case-cohort studies}
\runtitle{  Z-estimation and case-cohort design  }

\begin{aug}
\author{\fnms{Bin} \snm{Nan}\thanksref{t1}\ead[label=e1]{bnan@umich.edu}}
and
\author{\fnms{Jon A.} \snm{Wellner}\thanksref{t2}\ead[label=e2]{jaw@stat.washington.edu}}

\thankstext{t1}{Supported in part by NSF grant DMS-1007590 and NIH grant R01-AG036802.}
\thankstext{t2}{Supported in part by NSF grant DMS--1104832, NI-AID grant 2R01 AI291968-04,
     and the Alexander von Humboldt Foundation.}
%\thankstext{t3}{Second supporter of the project}
\runauthor{B. Nan and J. Wellner}

  \affiliation{University of Michigan and University of Washington}

\address{Department of Biostatistics \\ University of Michigan \\
Ann Arbor, MI 48109-2029\\
\printead{e1}\\
\phantom{E-mail:\ }}

\address{Department of Statistics\\
University of Washington\\
Seattle, WA 98195-4322\\
\printead{e2}\\
\phantom{E-mail:\ }}

\end{aug}

\begin{abstract}
Case-cohort design, an outcome-dependent sampling design for censored survival data, is increasingly used in biomedical research. The development of asymptotic theory for a case-cohort design in the current literature primarily relies on counting process stochastic integrals. Such an approach, however, is rather limited and lacks theoretical justification for outcome-dependent weighted methods due to non-predictability. Instead of stochastic integrals, we derive asymptotic properties for case-cohort studies based on a general Z-estimation theory for semiparametric models with bundled parameters using modern empirical processes. Both the Cox model and the additive hazards model with time-dependent covariates are considered.

\end{abstract}

%\begin{keyword}[class=AMS]
%\kwd[Primary ]{???}
%\kwd{???}
%\kwd[; secondary ]{???}
%\end{keyword}


\end{frontmatter}

\newpage

\section{Introduction}

Case-cohort designs, originally proposed by \cite{prentice:86} for right-censored survival data, are very useful in large epidemiologic
cohort studies, and their applications are increasingly common in
biomedical research. In a case-cohort study, complete data are
only obtained for all failures observed during follow-up and for a
sub-sample, called the subcohort, of the entire cohort. The subcohort can be a simple random or stratified sub-sample. Such a design is cost-effective for studies of rare events, and has been extended to other models including the additive hazards model \citep{kulich-lin:00}, transformation models \citep{chen-zucker:09, kong-etal:04, lu-tsiatis:06}, and the accelerated failure time model \citep{nan-kalbfleisch-yu:09, nan-yu-kalbfleisch:06}, and also to other censoring mechanisms \citep{li-etal:08, li-nan:11}, among many others.

For right-censored data, the pseudo likelihood approach of \cite{self-prentice:88} constructs risk sets from subcohort only, thus the counting process martingale theory is naturally applicable for deriving the asymptotic properties for the Cox-type regression models. This same strategy can be applied to some other regression models for right-censored data, for example, the accelerated failure time model studied by \cite{nan-yu-kalbfleisch:06}. Since complete information is also observed for all the failures, constructing risk sets from all observed data including failures outside the subcohort would yield more efficient estimation. This has been observed by many authors, for example, \cite{borgan-etal:00, chen-lo:99, chen-zucker:09, kalbfleisch-lawless:88, kulich-lin:00, kulich-lin:04, nan-kalbfleisch-yu:09}. The development of corresponding asymptotic theories has been primarily based on calculations of counting process stochastic integrals. Such a method, however, lacks theoretical justification because the integrands of those stochastic integrals are not predicable, not even adapted with respect to any filtration generated from the history.

To overcome this technical hurdle, we consider a general semiparametric Z-estimation method for {\sl bundled parameters} using empirical process theory, see e.g. \cite{vaart-wellner:96, vaart-wellner:07}. Our approach does not
use the stochastic integral formulation, thus there is no predictability requirement.
The main body of the article is as follows. In Section 2, we introduce a general asymptotic theory for semiparametric Z-estimation with bundled parameters. We then apply the Z-estimation theory to prove the asymptotic
properties for case-cohort studies in Section 3. Both the Cox model and the additive hazards model with time-dependent covariates will be considered. We
make some concluding remarks in Section 4.

\section{Semiparametric Z-estimation for bundled parameters}

Let $\theta \in \Theta \subset \mathbb{R}^d$ be the parameter of
interest, and $\eta: {\cal X}\times \Theta \to \mathbb{R}^J$ be
infinite dimensional nuisance parameter(s) in a Banach space ${\cal
H} \equiv \{ (x,\theta) \mapsto \eta(x,\theta) \in \mathbb{R}^J: x \in {\cal X}, \theta \in \Theta\}$. Such a parametrization allows the nuisance parameter to be a function of the parameter of interest, thus the two types of parameters are {\sl bundled together}, a terminology originally used by \cite{huang-wellner:97} and further studied by, for example, \cite{ding-nan:11}. Denote the random map
${\cal X}^n \mapsto \mathbb{R}^d$ with $n$
observations $X_1, \dots, X_n$ as
\begin{eqnarray}
\Psi_n(\theta; \eta)  \ \equiv \ \Psi_n(X_1, \dots, X_n; \theta, \eta(\cdot;\theta)) \, , \label{GeneralPsi-n}
\end{eqnarray}
which becomes an estimating
function for $\theta$ when $\eta$ is given or replaced by its
estimator. For independent and identically distributed (i.i.d.) observations $X_1, \dots, X_n$, very often $\Psi_n(\theta, \eta)$ takes the following form:
\begin{eqnarray}\Psi_n(\theta, \eta)  \ = \
\frac{1}{n}\sum_{i=1}^n \psi(X_i; \theta, \eta(\cdot;\theta)),  \label{LinearPsi-n}
\end{eqnarray}
where $\psi(\theta, \eta)  \ \equiv \ \psi(X; \theta, \eta(\cdot;\theta))$ is a random map ${\cal X} \mapsto \mathbb{R}^d$ with a single observation $X$.

Here we use the term ``nuisance parameter" in a
rather loose sense. It does not need to be an actual parameter (for example,
the baseline hazard function in the Cox model) in the original
parametrization of the distribution of $X$. Broadly speaking, it is
an unknown quantity in the estimating function in addition to the
parameter of interest. The unknown quantity $\eta$ as a function of $\theta$ needs to
be estimated prior to estimating $\theta$. We call the solution to
$\Psi_n(\theta, \hat{\eta}_n(\cdot;\theta))=0$ the
Z-estimator for $\theta$, where $\hat{\eta}_n$ is some estimator for
$\eta$. This type of generalization has been considered in the econometrics literature; see for example, \cite{newey:94, chen-etal:03}. We provide slightly modified results of \cite{chen-etal:03} with a focus on Z-estimation in the following lemmas, which we will use for the estimates in case-cohort studies we
consider in this article. Proofs of the lemmas are provided in the Appendix.

Let $\theta_0$ denote the true value of $\theta$ and $\eta_0$ be the
true functional form of $\eta$. Let $\Psi(\theta, \eta)$ be a
deterministic function, which usually denotes the limit of
$\Psi_n(\theta, \eta)$ as $n \to \infty$. We use $p^*$ to denote
``in outer probability", and refer its definition and detailed
discussion to \cite{vaart-wellner:96}.  Note
that all the lemmas in this section do not require i.i.d. data,
though data in the case-cohort studies we consider are assumed to be i.i.d. Let $| \cdot |$ be the Euclidian norm.
Let $\|
\cdot \|$ be the supremum of a norm or semi-norm taking over all $\theta \in \Theta$, that is $\|\eta\| = \sup_{\theta \in \Theta} \rho(\eta(\cdot;\theta))$ for some norm or semi-norm $\rho$; for example, $\rho(\eta(\cdot;\theta)) = \sup_{x \in {\cal X}} |\eta(x;\theta)|$, which gives $\|\eta\| = \sup_{\theta \in \Theta} \sup_{x \in {\cal X}} |\eta(x;\theta)|$.

{\lem \label{Lemma 2.1} (Consistency.) Suppose
$\theta_0$ is the unique solution to  $\Psi(\theta,
\eta_0(\cdot;\theta)) = 0$ in the parameter space $\Theta$ and
$\hat{\eta}_n$ is an estimator of $\eta_0$ such that $\|
\hat{\eta}_n - \eta_0 \| = o_{p^*}(1)$. If
\begin{eqnarray}
\sup_{\theta \in \Theta, \| \eta - \eta_0 \| \leq \delta_n}
\frac{|\Psi_n(\theta, \eta(\cdot;\theta)) - \Psi(\theta,
\eta_0(\cdot;\theta)) |}{1 + |\Psi_n(\theta, \eta(\cdot;\theta))| +
|\Psi(\theta, \eta_0(\cdot;\theta))|} = o_{p^*}(1)
\label{ConsistencyCondition}
\end{eqnarray}
for every sequence $\{\delta_n\} \downarrow 0$, then
$\hat{\theta}_n$ satisfying $\Psi_n(\hat{\theta}_n,
\hat{\eta}_n(\cdot;\hat{\theta}_n)) = o_{p^*}(1)$ converges in outer
probability to $\theta_0$. }

\medskip

Since consistency is a global property, so our main condition,
equation (\ref{ConsistencyCondition}), is therefore necessarily
global, that is, the supremum is taken over all of $\Theta$. The $p^{*}$
in equation (\ref{ConsistencyCondition}) indicates that the
left-hand side converges to $0$ in outer probability in case that
the term on the left is not Borel measurable. It is a stronger
condition to require that the convergence holds when the denominator
is replaced by 1. The purpose of adding an extra term in the
denominator is to control the numerator when it blows up to infinity
for some $\theta\in\Theta$.

%The above equation (\ref{ConsistencyCondition}) is more general than
%Assumption 5.4 in \cite{newey:94}, the condition for consistency,
%where $\psi(X; \theta, \eta_0)$ (not the deterministic function
%$\Psi$) is required to be continuous in $\theta$ with an envelope
%function $b(X)$, and for all $\theta \in \Theta$, $\|\psi(X; \theta,
%\eta) - \psi(X; \theta, \eta_0)\| \leq \tilde{b}(X) \|\eta -
%\eta_0\|^\epsilon$ for some norm $\|\cdot\|$ and $\epsilon>0$. The
%additive hazards model studied in Section 3.2 satisfies equation
%(\ref{ConsistencyCondition}), but clearly violates Newey's
%Assumption 5.4. Newey's Assumption 5.5 is the existence of a unique
%solution to $\Psi(\theta, \eta_0)=0$, which is also a sufficient
%condition in Lemma \ref{Lemma 2.1}.

{\lem \label{Lemma 2.2} (Rate of convergence and asymptotic
representation.)  Let ${\cal H}_0 = \{\eta(x;\theta): x\in {\cal X},
\theta \in \Theta_0\}$ be a collection of functions that are
continuously differentiable in $\theta$ for all $x \in {\cal X}$
with bounded derivative matrices $\{\dot\eta(\cdot; \theta)\}$,
where $\Theta_0 \subset \Theta$ is a neighborhood of $\theta_0$.
Suppose that $\hat{\theta}_n$ satisfying $\Psi_n(\hat{\theta}_n,
\hat{\eta}_n(\cdot; \hat{\theta}_n)) = o_{p^*}(n^{-1/2})$ is a
consistent estimator of $\theta_0$ that is the unique solution to the equation 
$\Psi(\theta, \eta_0(\cdot;\theta)) = 0$ in $\Theta$, and that
$\hat{\eta}_n \in {\cal H}_0$ is an estimator of $\eta_0 \in {\cal
H}_0$ satisfying $\| \hat{\eta}_n - \eta_0 \| = O_{p^*}(n^{-\beta})$
for some $\beta>0$. Suppose the following four conditions are
satisfied:

\medskip
(i) (Stochastic equicontinuity.)
$$
\frac{|n^{1/2}(\Psi_n - \Psi)(\hat{\theta}_n, \hat{\eta}_n(\cdot;
\hat{\theta}_n)) - n^{1/2}(\Psi_n - \Psi)(\theta_0, \eta_0(\cdot;
\theta_0))|}{1 + n^{1/2}|\Psi_n(\hat{\theta}_n, \hat{\eta}_n(\cdot;
\hat{\theta}_n))| + n^{1/2}|\Psi(\hat{\theta}_n, \hat{\eta}_n(\cdot;
\hat{\theta}_n))|} = o_{p^*}(1) \ .
$$

\medskip
(ii) $n^{1/2}\Psi_n(\theta_0, \eta_0(\cdot; \theta_0)) =
O_{p^*}(1).$

\medskip
(iii) (Smoothness.) (a) If $\beta=1/2$, the function $\Psi(\theta, \eta(\cdot; \theta)):
\Theta_0 \times {\cal H}_0 \to \mathbb{R}^d$ is Fr\'{e}chet
differentiable at $(\theta_0, \eta_0(\cdot;\theta_0))$,
 i.e., there exists a continuous $d \times d$ matrix
$\dot{\Psi}_1(\theta_0, \eta_0(\cdot;\theta_0))$ and a continuous
linear functional $\dot{\Psi}_2(\theta_0, \eta_0(\cdot;\theta_0))$
such that
\begin{eqnarray}
&& |\Psi(\theta, \eta(\cdot;\theta)) -
\Psi(\theta_0, \eta_0(\cdot;\theta_0)) \nonumber \\
&& \qquad \qquad - \
\{\dot{\Psi}_1(\theta_0,\eta_0(\cdot;\theta_0))+\dot{\Psi}_2(\theta_0,
\eta_0(\cdot;\theta_0))[\dot{\eta}_0(\cdot;\theta_0)]\}(\theta -
\theta_0)
\nonumber \\
&&  \qquad \qquad  - \
 \dot{\Psi}_2(\theta_0, \eta_0(\cdot;\theta_0))[(\eta - \eta_0)(\cdot;\theta_0)]| \nonumber \\
 && \qquad  = \ o(|\theta -
\theta_0|) + o(\|\eta - \eta_0\|) \, ; \label{hu:smooth}
\end{eqnarray}
or (b) if $0 < \beta < 1/2$, for some $\alpha > 1$ satisfying $\alpha \beta > 1/2$ we have
\begin{eqnarray}
&& |\Psi(\theta, \eta(\cdot;\theta)) -
\Psi(\theta_0, \eta_0(\cdot;\theta_0)) \nonumber \\
&& \qquad \qquad - \
\{\dot{\Psi}_1(\theta_0,\eta_0(\cdot;\theta_0))+\dot{\Psi}_2(\theta_0,
\eta_0(\cdot;\theta_0))[\dot{\eta}_0(\cdot;\theta_0)]\}(\theta -
\theta_0)
\nonumber \\
&& \qquad \qquad - \
 \dot{\Psi}_2(\theta_0, \eta_0(\cdot;\theta_0))[(\eta -
 \eta_0)(\cdot;\theta_0)]| \nonumber \\
 &&  \qquad = \ o(|\theta - \theta_0|) + O(\|\eta -
\eta_0\|^\alpha) \, . \label{huang:smooth}
\end{eqnarray}
Here the subscripts 1 and 2 correspond to the first and the second
arguments in $\Psi(\cdot,\cdot)$, respectively, and we assume that the matrix
$$A = - \dot{\Psi}_1(\theta_0,\eta_0(\cdot;\theta_0)) - \dot{\Psi}_2(\theta_0,
\eta_0(\cdot;\theta_0))[\dot{\eta}_0(\cdot;\theta_0)]$$ is
nonsingular.

\medskip
(iv) $n^{1/2}\dot{\Psi}_2(\theta_0,
\eta_0(\cdot;\theta_0))[(\hat{\eta}_n - \eta_0)(\cdot;\theta_0)] =
O_{p^*}(1)$.

\medskip
Then $\hat{\theta}_n$ is $n^{1/2}$-consistent and further we have
\begin{eqnarray}
n^{1/2}(\hat{\theta}_n - \theta_0) &=&   A^{-1} n^{1/2}\Big\{ (\Psi_n - \Psi)(\theta_0,
\eta_0(\cdot;\theta_0)) \nonumber \\
&& \quad  + \ \dot{\Psi}_2(\theta_0,
\eta_0(\cdot;\theta_0))[(\hat{\eta}_n - \eta_0)(\cdot;\theta_0)]
\Big\} +  o_{p^*}(1) .\label{AsymptoticRepresentation}
\end{eqnarray}
}

{\sc Remark:} For i.i.d. data, Condition (i) in Lemma \ref{Lemma 2.2} holds if
the class of functions $\{\psi(\theta, \eta): |\theta - \theta_0| <
\delta, \|\eta - \eta_0 \| < \delta \}$ is Donsker for some $\delta
> 0$ and satisfies $E_0|\psi(\theta, \eta; X) - \psi(\theta_0, \eta_0
; X)|^2 \to 0$ as $|\theta - \theta_0| \to 0$ and $\|\eta - \eta_0\|
\to 0$ (see e.g. Corollary 2.3.12 of \cite{vaart-wellner:96}, page
115). Though simpler, this is stronger than Condition (i). Condition
(ii) holds automatically for i.i.d. data if $E_0|\psi(\theta_0,
\eta_0)|^2 < \infty$ and $\Psi_n$ takes the form in (\ref{LinearPsi-n}). In Condition (iii),
$\{\dot{\Psi}_1(\theta_0,\eta_0(\cdot;\theta_0))+\dot{\Psi}_2(\theta_0,
\eta_0(\cdot;\theta_0))[\dot{\eta}_0(\cdot;\theta_0)]\}(\theta -
\theta_0)$ is obtained by the chain rule, which is the usual inner
product of a $d\times d$ matrix and a $d\times 1$ vector; whereas
$\dot{\Psi}_2(\theta_{0},\eta_{0}(\cdot;\theta_0))[(\eta-\eta_{0})(\cdot;\theta_0)]
=\sum_{j=1}^{J}\dot{\Psi}_{2_{j}}(\theta_{0},\eta_{0}(\cdot;\theta_0))[(\eta_{j}-\eta_{0j})(\cdot;\theta_0)]$,
here $J$ is the number of infinite dimensional parameters contained
in $\eta$, is the sum of separate terms with each
$\dot{\Psi}_{2_{j}}$ being a bounded linear functional that brings
$\eta-\eta_{0}$ to a real number, where $\eta$ is close to
$\eta_{0}$ in $n^{\beta}$-rate for some $\beta>0$. Note that
equation (\ref{huang:smooth}) is indeed a stronger condition than
equation (\ref{hu:smooth}). Proposition 1 of \cite{bkrw:93}, page
455, provides useful tools for checking Fr\'{e}chet
differentiability for infinite-dimensional parameters. Condition
(iv) holds automatically under (iii) if $\hat{\eta}_n$ is
$n^{1/2}$-consistent, but may require extensive work for slower than
root-$n$ convergence rate, see e.g. \cite{wong-severini:91} and \cite{huang-wellner:95}.
In view of the structure of equation
(\ref{AsymptoticRepresentation}), the asymptotic distribution of
$n^{1/2}(\hat{\theta}_{n}-\theta_{0})$ is determined by the
asymptotic joint distribution of the random variables
$n^{1/2}(\Psi_{n}-\Psi)(\theta_{0}, \eta_{0}(\cdot;\theta_0))$ and
$n^{1/2}\dot{\Psi}_2(\theta_{0},
\eta_{0}(\cdot;\theta_0))[(\hat{\eta}_{n}-\eta_{0})(\cdot;\theta_0)]$,
particularly if the asymptotic joint distribution is multivariate
Gaussian.

%In contrast to the above Condition (i) for the estimating function
%$\Psi_n$, \cite{newey:94} requires stochastic equicontinuity for
%the derivative of the estimating function to $\eta$ in his
%Assumption 5.2. The above Condition (i) is a more natural assumption
%and in general should be easier to verify than Newey's Assumption
%5.2. Newey's Assumptions 5.1 and 5.6(ii) are for the
%differentiability of the random map $\psi(\theta, \eta)$, while our
%Condition (iii) requires differentiability of the deterministic
%function $\Psi(\theta, \eta)$. So our method can handle nonsmooth
%estimating functions as long as $\Psi$ is smooth. Newey's Assumption
%5.3 together with interchangeability of differentiation and
%integration implies that $n^{1/2}\dot{\Psi}_2(\theta_0,
%\eta_0)(\hat{\eta}_n - \eta_0)$ is asymptotically linear, which is
%stronger than our Condition (iv), but something we would expect to
%hold in order to easily obtain asymptotic normality of
%$n^{1/2}(\hat{\theta}_n - \theta_0)$ from the above equation
%(\ref{AsymptoticRepresentation}). His Assumption 5.6(iv) renders our
%Condition (ii) for i.i.d. data as discussed in the above Remark, and
%Assumption 5.6(iii) is also assumed in our Condition (iii).

\medskip

In the case that $\eta$ is free of $\theta$, we have $\dot{\eta} =
0$. Then Lemma \ref{Lemma 2.2} reduces to the following corollary that was
studied by \cite{hu:98}. The corollary is particularly useful for the case-cohort additive hazards model in the next section. Now we replace $\dot{\Psi}_1$ by
$\dot{\Psi}_\theta$ and $\dot{\Psi}_2$ by $\dot{\Psi}_\eta$ without
causing any confusion, and the notation $\| \cdot \|$ becomes a norm.

%\bigskip

{\cor \label{Corollary 2.3} (Rate of convergence and asymptotic
representation.) Suppose that $\hat{\theta}_n$ satisfying
$\Psi_n(\hat{\theta}_n, \hat{\eta}_n) = o_{p^*}(n^{-1/2})$ is a
consistent estimator of $\theta_0$ that is the unique solution to
$\Psi(\theta, \eta_0) = 0$ in $\Theta$, and that $\hat{\eta}_n$ is
an estimator of $\eta_0$ satisfying $\| \hat{\eta}_n - \eta_0 \| =
O_{p^*}(n^{-\beta})$ for some $\beta>0$. Suppose the following four
conditions are satisfied:

\medskip
(i) (Stochastic equicontinuity.)
$$
\frac{|n^{1/2}(\Psi_n - \Psi)(\hat{\theta}_n, \hat{\eta}_n) -
n^{1/2}(\Psi_n - \Psi)(\theta_0, \eta_0)|}{1 +
n^{1/2}|\Psi_n(\hat{\theta}_n, \hat{\eta}_n)| +
n^{1/2}|\Psi(\hat{\theta}_n, \hat{\eta}_n)|} = o_{p^*}(1) \ .
$$

\medskip
(ii) $n^{1/2}\Psi_n(\theta_0, \eta_0) = O_{p^*}(1).$

\medskip
(iii) (Smoothness.) (a) If $\beta=1/2$, function $\Psi(\theta, \eta)$ is Fr\'{e}chet
differentiable at $(\theta_0, \eta_0)$, i.e., there
exists a continuous and nonsingular $d \times d$ matrix
$\dot{\Psi}_\theta(\theta_0, \eta_0)$ and a continuous linear
functional $\dot{\Psi}_\eta(\theta_0, \eta_0)$ such that
\begin{eqnarray} \label{hu:smooth2}
&& |\Psi(\theta, \eta) - \Psi(\theta_0, \eta_0) -
\dot{\Psi}_\theta(\theta - \theta_0) - \dot{\Psi}_\eta(\theta_0,
\eta_0)[\eta - \eta_0]| \\
&& \qquad = \ o(|\theta - \theta_0|) + o(\|\eta -
\eta_0\|) ; \nonumber
\end{eqnarray}
or (b) if $0 < \beta < 1/2$, for some $\alpha > 1$ satisfying $\alpha \beta > 1/2$ we have
\begin{eqnarray} \label{huang:smooth2}
&& |\Psi(\theta, \eta) - \Psi(\theta_0, \eta_0) -
\dot{\Psi}_\theta(\theta - \theta_0) - \dot{\Psi}_\eta(\theta_0,
\eta_0)[\eta - \eta_0]| \\
&& \qquad = \ o(|\theta - \theta_0|) \! + \! O(\|\eta -
\eta_0\|^\alpha) . \nonumber
\end{eqnarray}

\medskip
(iv) $n^{1/2}\dot{\Psi}_\eta(\theta_0, \eta_0)[\hat{\eta}_n -
\eta_0] = O_{p^*}(1)$.

\medskip
Then $\hat{\theta}_n$ is $n^{1/2}$-consistent and further we have
\begin{eqnarray}
&& n^{1/2}(\hat{\theta}_n - \theta_0)  \label{AsymptoticRepresentation2} \\
&& \qquad =  \Big\{ - \dot{\Psi}_\theta(\theta_0,
\eta_0)\Big\}^{-1} n^{1/2}\Big\{ (\Psi_n - \Psi)(\theta_0, \eta_0) +
\dot{\Psi}_\eta(\theta_0, \eta_0)[\hat{\eta}_n - \eta_0] \Big\}  \nonumber \\
&& \qquad \qquad  + \
o_{p^*}(1).\nonumber
\end{eqnarray}

}

%\medskip

%For the examples  we consider in the
%following section, we need to carefully characterize the term
%$\dot{\Psi}_2(\theta_{0},
%\eta_{0}(\cdot;\theta_0))\left(\hat{\eta}_{n}-\eta_{0}\right)(\cdot;\theta_0)$
%in the asymptotic representation (\ref{AsymptoticRepresentation}),
%or $\dot{\Psi}_\eta(\theta_{0},
%\eta_{0})\left(\hat{\eta}_{n}-\eta_{0}\right)$ in the asymptotic
%representation (\ref{AsymptoticRepresentation2}). In the next
%section we show that this term can be asymptotically expressed as
%the average of i.i.d. zero mean random variables. Thus the
%asymptotic distribution of $n^{1/2}(\hat{\theta}_{n}-\theta_{0})$ is
%easily determined.

\section{Case-Cohort Studies}

We consider two models that are used for
analyzing case-cohort data: the Cox model and the
additive hazards model. Let $X$ be the generic random variable
that consists of several random variables. Let $T$ be the failure time
and $C$ the censoring time, we only observe $Y = {\rm min}(T,C)$
and the failure indicator $\Delta = 1(T \leq C)$. Let $Z(\cdot)$ be the
$d$-dimensional covariate process and $\bar Z(t)$ be the covariate history up to time $t$. We assume that for all $t$, events $\{T \ge t\}$ and $\{C \ge t\}$ are conditionally independent given $\bar Z(t)$, and both are independent of $\{\bar Z(s): s > t\}$. In other words, $Z(\cdot)$ is an external covariate, see \cite{kalbfleisch-prentice:02}. Suppose potentially we would have $n$ i.i.d.
copies of $(Y, \Delta, \bar Z(Y))$ in the full cohort, but we only observe
$\bar Z(Y)$ for all failures and subjects in the
subcohort that is a sub-sample of the entire cohort. The subcohort
may be selected using a variety of sampling schemes including the
simple random sampling and the stratified sampling based on some auxiliary variable $Z^*(\cdot)$ that can be a subset of $Z(\cdot)$, may or may not be time-dependent, and is available to everyone in the
cohort. We focus on the independent Bernoulli sampling method for
selecting the subcohort by which a coin is flipped for each subject $i$ in
the cohort with a given success probability $\pi_i$ that may depend
on $Z^*_i$. For finite population sampling methods, as applied in
\cite{breslow-wellner:06}, we expect the weighted bootstrap
empirical process theory of \cite{p-wellner:93} to be a useful tool
to verify conditions in Lemmas \ref{Lemma 2.1} and \ref{Lemma 2.2}. See \cite{saegusa-wellner:12} for a related problem using the weighted bootstrap empirical process theory.

Let $R_i$ be the subcohort indicator that equals 1 if
the $i$th subject is selected into the subcohort and 0 otherwise.
Then $\pi_i = P(R_i=1|Z^*_i)$. Thus the observed data in such a
case-cohort study are i.i.d. and the missing data mechanism is missing at random \citep{little-rubin:02}.
The following is a set of common regularity conditions for both models.

\medskip
{\sc Assumption} (A): The sample paths of $Z(\cdot) \in {\cal Z}$ are bounded with bounded variation, and the parameter space $\Theta$ is compact.

\smallskip
{\sc Assumption} (B): The conditional distribution of $T$ given $\bar Z(\cdot)$
possesses a continuous Lebesgue density.

\smallskip
{\sc Assumption} (C): The study stops at a finite time $\tau > 0$
such that, for constants $\sigma_1$ and $\sigma_2$, $\inf_{z \in {\cal Z}}P(C \geq \tau | \bar Z(\tau) = \bar z(\tau)) = \sigma_1
> 0$ and $\inf_{z \in {\cal Z}}P(T > \tau | \bar Z(\tau) = \bar z(\tau)) = \sigma_2 \in (0,1).$

\smallskip
{\sc Assumption} (D): The map $\Psi(\theta, \eta(\cdot;\theta)) =
P\psi(\theta, \eta(\cdot;\theta))$ is Fr\'{e}chet differentiable at
$(\theta_0, \eta_0(\cdot;\theta_0))$ with a nonsingular partial
derivative with respect to $\theta$ at $(\theta_0, \eta_0(\cdot;\theta_0))$.

\smallskip

{\sc Assumption} (E): In case-cohort studies, data are missing at
random with $\pi_i \geq \sigma_3
> 0$ for all $i$ and a constant $\sigma_3$.

\medskip

Note that the assumption of compact $\Theta$ is only for technical convenience, which is unnecessarily strong. Later we will see that for the additive hazards model, $\eta$ is free of $\theta$. The following is some standard empirical
process notation that we will use in the rest of the paper. Suppose
$X_1, \dots, X_n$ are i.i.d. $p$-dimensional random variables that
follow the distribution $P$ on a measurable space $({\cal X}, {\cal
A})$. For a measurable function $ f: {\cal X} \mapsto \mathbb{R} $,
we denote
\begin{eqnarray*}
\mathbb{P}_n f &=& \frac{1}{n}\sum_{i=1}^n f(X_i) \, , \quad
P f  \ = \ \int f dP \, , \quad {\rm and} \\
\mathbb{G}_n f &=& n^{-1/2}\sum_{i=1}^n \{ f(X_i) - Pf \} \ = \
n^{1/2}(\mathbb{P}_n - P)f \, .
\end{eqnarray*}
%Hence the notation $E$ and $P$ are exchangeable when they denote
%expectations.
Function $f$ can be replaced by a random function $x
\mapsto \hat{f}_n(x; X_1, \dots, X_n)$. Thus,
$$
\mathbb{P}_n \hat{f}_n = \frac{1}{n}\sum_{i=1}^n \hat{f}_n(X_i;
X_1, \dots,X_n) \, , \quad  P \hat{f}_n  \ = \ \int \hat{f}_n(x;
X_1, \dots,X_n) dP(x) \, ,$$  $$ \mbox{and} \quad \mathbb{G}_n \hat{f}_n
= n^{-1/2}\sum_{i=1}^n \{ \hat{f}_n(X_i; X_1, \dots, X_n ) - P
\hat{f}_n \} \ = \ n^{1/2}(\mathbb{P}_n - P) \hat{f}_n \, .$$

\subsection{Case-cohort study: the Cox model}

For the Cox model with external time-dependent covariates, we have
$$
\lambda(t|\bar Z(t)) = \lambda_0(t) e^{\theta_0' Z(t)}
$$
and
$$
1 - F_{T|\bar Z(\tau)}(t|\bar z(t)) = \exp\left\{-\int_0^t e^{\theta'_0 z(s)}d\Lambda_0(s)
\right\} \, ,
$$
where $F_{T|\bar Z(\tau)}$ is the conditional distribution function of $T$
given $\bar Z(\tau)$, $\Lambda_0$ is the baseline cumulative hazard function,
and $\theta_0$ is the parameter of interest. We define the following random map
\begin{eqnarray}
\Psi_n(\theta, \eta) = \frac{1}{n}\sum_{i=1}^n \Omega_i \{Z_i(Y_i) -
\eta(Y_i;\theta)\}\Delta_i \, , \label{Cox}
\end{eqnarray}
with true $\eta$ given by
$$
\eta_0(t; \theta) = \frac{E\{Z(t)e^{\theta'Z(t)}1(Y \geq
t)\}}{E\{e^{\theta'Z(t)}1(Y \geq t)\}} \, ,
$$
where $\Omega_i$ are diagonal weight matrices with subject and
covariate specific random weights on the diag that have expectation
1 given complete data $X_i = (Y_i, \Delta_i, \bar Z_i(Y_i), \bar Z_i^*(Y_i))$. By
choosing a weight matrix, we are allowed to weight each component of
$\psi(X_i; \theta, \eta)$ differently, as in \cite{kulich-lin:04}. For notational simplicity, we consider a scalar weight $\Omega_i$ in the rest of the article. The proofs for a matrix $\Omega_i$ are almost identical.
It has been shown by \cite{andersen-gill:82} that $E\psi(\theta_0, \eta_0(\cdot; \theta_0)) = E[ \{Z(Y)
- \eta_0(Y;\theta_0)\}\Delta] = 0$. The explicit functional form of
$\eta_0$ is unknown and needs to be estimated first in order to
estimate $\theta$ from (\ref{Cox}).

For full-cohort data,
$\Omega_i = 1$, and the partial
likelihood estimating function is
\begin{eqnarray}
\Psi_n(\theta, \hat\eta_n^F) = \frac{1}{n} \sum_{i=1}^n \{ Z_i(Y_i) -
\hat\eta_n^F (Y_i; \theta) \} \Delta_i \, , \label{FullCox}
\end{eqnarray}
where $\hat\eta_n^F$ is an estimator of $\eta_0$ using full data, which has the following
form:
$$
\hat\eta_n^F (t; \theta) = \frac{\sum_{j=1}^n  Z_j(t) e^{\theta'Z_j(t)}1(Y_j
\geq t)}{\sum_{j=1}^n  e^{\theta'Z_j(t)}1(Y_j \geq t)} \, .
$$

For case-cohort data where the subcohort is a sub-sample of the
entire cohort selected with a constant probability $\pi_i$ for all
$i$, also with $\Omega_i = 1$, the pseudo-likelihood
estimating function of \cite{self-prentice:88} is
\begin{eqnarray}
\Psi_n(\theta, \hat\eta_n^{SP}) = \frac{1}{n} \sum_{i=1}^n \{ Z_i(Y_i) -
\hat\eta_n^{SP}(Y_i; \theta) \} \Delta_i \, , \label{CCCox}
\end{eqnarray}
where $\hat\eta_n^{SP}$ is an estimator of $\eta_0$ considered by \cite{self-prentice:88} using the
subcohort data only, which has the following form:
$$
\hat\eta_n^{SP}(t; \theta) = \frac{\sum_{j\in {\cal SC}} Z_j(t)
e^{\theta'Z_j(t)}1(Y_j \geq t)}{\sum_{j \in {\cal SC}}
e^{\theta'Z_j(t)}1(Y_j \geq t)} \, .
$$
Here ${\cal SC}$ denotes the set of subjects in the subcohort.

In order to improve efficiency, the subcohort can be chosen by
stratified sampling, and furthermore, it is tempting to include
failures outside the subcohort to estimate $\eta_0$, see e.g.
\cite{kalbfleisch-lawless:88}. The corresponding estimating
function then becomes
\begin{eqnarray}
\Psi_n(\theta, \hat\eta_n^{W}) = \frac{1}{n} \sum_{i=1}^n \Omega_i(Y_i) \{
Z_i(Y_i) - \hat\eta_n^{W}(Y_i; \theta) \} \Delta_i \, , \label{WCCCox}
\end{eqnarray}
where $\hat\eta_n^{W}$ is a weighted estimator of $\eta_0$ with
 the following form
$$
\hat\eta_n^{W}(t; \theta) = \frac{\sum_{j = 1}^n W_j(t)
Z_j(t) e^{\theta'Z_j(t)}1(Y_j \geq t)}{\sum_{j=1}^n W_j(t)
e^{\theta'Z_j(t)}1(Y_j \geq t)}   \, .
$$
Here $W_i$ could also be diagonal weight matrices with subject and
covariate specific random weights on the diag. Again for notational simplicity, we consider scalar $W_i$, which may or may not
equal to $\Omega_i$. We also require that $W_i$
have expectation 1 given complete data $X_i = (Y_i, \Delta_i, Z_i(\cdot),
Z_i^*(\cdot))$. We consider a broad class of weighted problems by allowing both weights $\Omega$ and $W$ to be time-dependent. The commonly used weights, originally proposed by
\cite{kalbfleisch-lawless:88}, are the inverse-probability weights
\begin{eqnarray}
W_i = \Delta_i + \frac{R_i}{\pi_i}(1-\Delta_i) \ ,
\label{Weights}
\end{eqnarray}
where $\pi_i$ can be time-dependent, see \cite{kulich-lin:04} for example.

Note that the estimating functions in (\ref{FullCox}) and
(\ref{CCCox}) can be expressed by using counting
process stochastic integrals and martingale theory applies in deriving
asymptotic properties of corresponding estimators, see e.g.
\cite{andersen-gill:82} and \cite{self-prentice:88}. Using a similar stochastic integral for the estimating
function (\ref{WCCCox}) with weights (\ref{Weights}), however, creates a
measurability problem because the integrand is no longer adapted to any
meaningful filtration (and hence not predictable).
See e.g. \cite{chung-williams:90} and \cite{protter:04} for detailed discussions on
stochastic integration.
In this article, instead of using stochastic integrals, we give a
rigorous proof of asymptotic properties of the estimators obtained
from the estimating function (\ref{WCCCox}) using the general
Z-estimation theory provided in Section 2.

It grants great flexibility in estimating $\theta$ from equation
(\ref{WCCCox}) to use two possibly different weights
$\Omega_i$ and $W_i$ . When $\Omega_i = W_i = 1$, the estimating
function $\Psi_n(\theta, \hat{\eta}_n^W(\cdot;\theta))$ reduces to
(\ref{FullCox}); that is, the partial likelihood estimating function
of \cite{cox:72} for full-cohort data. When $\Omega_i = 1$ and
$W_i = R_i/\pi_i$ with constant $\pi_i = \pi > 0$ for all
$i$, $\Psi_n(\theta, \hat{\eta}_n^W(\cdot;\theta))$ becomes (\ref{CCCox});
that is, the pseudo-likelihood estimating function of
\cite{self-prentice:88}. When $\Omega_i = W_i$ and they take the
form in (\ref{Weights}), $\Psi_n(\theta, \hat{\eta}_n^W(\cdot;\theta))$ is
equivalent to the weighted estimating function of
\cite{kalbfleisch-lawless:88}. When $\Omega_i = W_i = R_i^* /
\pi_i^*$, here $R_i^*$ is the indicator that equals 1 if
subject $i$ has complete data and 0 otherwise, and $\pi_i^* =
P(R_i^* =1 | X_i)$, $\Psi_n(\theta, \hat{\eta}_n^W(\cdot;\theta))$ becomes
the estimating function proposed by \cite{pugh:92}, which can be derived from a weighted likelihood method
for a two-phase design. The corresponding asymptotic properties
have been studied by \cite{breslow-wellner:06} for both
independent stratified Bernoulli sampling and finite population
stratified sampling when covariates are time-independent. To improve efficiency, \cite{kulich-lin:04} considered the estimating function $\Psi_n(\theta,
\hat{\eta}_n^W(\cdot;\theta))$ with $\Omega_i = 1$ and $W_i$ being time-dependent weights.
A clear advantage of introducing weights $\Omega_i$ in $\Psi_n(\theta,
\hat{\eta}_n^W(\cdot;\theta))$ is that it allows one to estimate $\theta$
from a data set in which some failures may have missing data, e.g.
the two-phase design studied by \cite{breslow-wellner:06}. This is
more general than a traditional case-cohort study which requires all failures
to be completely observed. It is obvious that all the above weights
are nonnegative and bounded, have unit conditional expectation given
complete data by Assumption (E), and are equal to zero if
corresponding covariates are missing. We will assume this holds
throughout the rest of the paper.

{\prop \label{Proposition 3.1}  Let $\hat{\eta}_n(t;\theta) =
\hat\eta_n^{W}(t;\theta)$ as in equation (\ref{WCCCox}). Suppose the weight process $W(t)$ has bounded sample paths of bounded variation. Then both $\hat\eta_n(t; \theta)$ and $\eta_0(t;\theta)$ belong to a Donsker class, and further we have
$\|\hat{\eta}_n - \eta_0 \| = O_{p^*}(n^{-1/2})$. }

\smallskip
{\sc Proof:} We consider one nuisance parameter $\eta$ for
simplicity. The vector $\eta$ can be dealt with by examining each of
its components. Define
\begin{eqnarray*}
D^{(0)}_n (t, \theta) & \equiv  & \mathbb{P}_n \Big\{W(t) e^{\theta'Z(t)}1(Y
\geq t) \Big\} , \\
d^{(0)} (t, \theta) & \equiv & P \Big\{W(t)
 e^{\theta'Z(t)}1(Y \geq t) \Big\}  = P \Big\{
 e^{\theta'Z(t)}1(Y \geq t) \Big\};
\end{eqnarray*}
and
\begin{eqnarray*}
D^{(1)}_n (t, \theta) &\!\!\! \equiv \!\!\! & \mathbb{P}_n \Big\{W(t) Z(t)e^{\theta'Z(t)} 1(Y
\geq t)\Big\} , \\
 d^{(1)} (t, \theta) & \equiv & P \Big\{W(t)  Z(t)
e^{\theta'Z(t)} 1(Y \geq t)\Big\} = P \Big\{ Z(t)
e^{\theta'Z(t)} 1(Y \geq t)\Big\}.
\end{eqnarray*}
Then we have
$$
\hat{\eta}_n(t; \theta) = \frac{D^{(1)}_n (t, \theta)}{D^{(0)}_n (t,
\theta)} \, , \quad \eta_0(t; \theta) = \frac{d^{(1)} (t,
\theta)}{d^{(0)} (t, \theta)} \, .
$$

Apparently the sets of functions ${\cal F}_0 = \{W(t) 1(Y \geq t)
e^{\theta'Z(t)}: 0 \leq t\leq \tau, \theta\in\Theta\}$ and ${\cal F}_1
= \{W(t) 1(Y \geq t) Z(t) e^{\theta'Z(t)}: 0 \leq t\leq \tau,
\theta\in\Theta\}$ are well-behaved and belong to Donsker classes,
see e.g. \cite{vaart-wellner:96}, Section 2.10. Hence we have that
$n^{1/2}\{D^{(k)}_n (t, \theta) - d^{(k)} (t, \theta)\}$ converge
weakly to zero mean Gaussian processes, and $\|D_n^{(k)} - d^{(k)}\|
= O_{p^*}(n^{-1/2})$, $k=0, 1$. Let $\bar{\cal F}_k$ be the closure
of ${\cal F}_k$, $k = 0, 1$, respectively, in which the convergence
is both pointwise and in $L_2(P)$. Then $D^{(k)}_n (t, \theta)$ and
$d^{(k)} (t, \theta)$ are in the convex hull of $\bar{\cal F}_k$, $k
= 0, 1$, and thus Donsker. See e.g. \cite{vaart-wellner:96},
Theorems 2.10.2 and 2.10.3. Hence both $\{\hat{\eta}_n(t; \theta)\}$
and $\{ \eta_0(t;\theta)\}$ are Donsker by
\cite{vaart-wellner:96}, Example 2.10.9, where $D_n^{(0)}$
and $d^{(0)}$ are bounded away (almost surely) from zero by Assumption (C).

Now we verify that $\hat{\eta}_n$ is $n^{1/2}$-consistent by the
following calculation:
\begin{eqnarray*}
&& n^{1/2} \{ \hat{\eta}_n(t;\theta) - \eta_0(t;\theta) \} \\
&& \qquad = \
n^{1/2} \Bigg[
\frac{1}{d^{(0)}(t,\theta)} \{D^{(1)}_n (t, \theta) -  d^{(1)} (t, \theta)\}  \\
&& \qquad \qquad \qquad - \ \frac{D^{(1)}_n(t,\theta)}{D^{(0)}_n(t,\theta)
d^{(0)}(t, \theta)}\{D^{(0)}_n (t, \theta) - d^{(0)} (t, \theta)\}
\Bigg] \\
&& \qquad = \ n^{1/2} \Bigg[
\frac{1}{d^{(0)}(t,\theta)} \{D^{(1)}_n (t, \theta) -  d^{(1)} (t, \theta)\}  \\
&& \qquad \qquad \qquad - \ \frac{d^{(1)}(t,\theta)}{ d^{(0)}(t,
\theta)^2}\{D^{(0)}_n (t, \theta) - d^{(0)} (t, \theta)\} \Bigg]
+ o_{p^*}(1) \\
&& \qquad = \ d^{(0)}(t, \theta)^{-1} \mathbb{G}_n \Big[ W(t)\{Z(t)-\eta_0(t;
\theta)\} e^{\theta'Z(t)}1(Y \geq t)\Big] + o_{p^*}(1) \, .
\end{eqnarray*}
Since the classes of functions $\{W(t)\}$, $\{1(Y\geq t)\}$, $\{Z(t)\}$,
and $\{e^{\theta'Z(t)}\}$ are all Donsker, and $\eta_0$ is a bounded
deterministic function, we know that the class $\{W(t)\{Z(t)-\eta_0(t; \theta)\}
e^{\theta'Z(t)}1(Y \geq t)\}$ is Donsker (see e.g.
\cite{vaart-wellner:96}, Section 2.10). We then obtain the desired
result. \hspace{5mm} $\Box$

{\prop \label{Proposition 3.2} Assume the conditions in  Proposition 3.1 and suppose the weight process $\Omega(t)$ also has bounded sample paths of bounded variation. Then the root of function (\ref{WCCCox})
denoted as $\hat{\theta}_n$ is a consistent estimator of $\theta_0$.
}

\smallskip

{\sc Proof:} We prove by verifying conditions in Lemma \ref{Lemma 2.1}. The
uniqueness of $\theta_0$ as a root of $\Psi(\theta,
\eta_0(\cdot;\theta))$ is proved by \cite{andersen-gill:82}, here
$\Psi(\theta, \eta_0(\cdot;\theta))$ corresponds to the derivative
of the limit of their function (2.7). The uniform consistency of
$\hat{\eta}_n$ is given by Proposition 3.1. Now we verify condition
(\ref{ConsistencyCondition}) by the following argument. Again we
consider one-dimensional $\theta$ for simplicity. Suppose that $\Omega_i <
K < \infty$ for all $i$ for a constant $K$. Let $\|\eta - \eta_0 \|
\leq \delta_n \downarrow 0$. Then we have
\begin{eqnarray*}
&& \hspace{-0.2in}|\Psi_n(\theta, \eta(\cdot;\theta)) - \Psi(\theta,
\eta_0(\cdot;\theta))| \\
&& \qquad  = \ \Big|\mathbb{P}_n [ \Omega(Y) \{Z(Y) - \eta(Y;
\theta)\}\Delta] - P[ \Omega(Y) \{Z(Y)
- \eta_0(Y; \theta)\}\Delta] \Big| \\
&& \qquad  \leq \ \Big| \mathbb{P}_n [\Omega(Y) Z(Y) \Delta] - P[\Omega(Y)
 Z(Y)\Delta] \Big|   \\
 && \qquad \qquad \qquad + \  \Big|\mathbb{P}_n \left[ \Omega(Y) \{ \eta(Y;
\theta) - \eta_0(Y;
\theta)\} \Delta \right] \Big| \\
&& \qquad \qquad \qquad  + \ \Big| (\mathbb{P}_n - P) [\Omega(Y)
\eta_0(Y; \theta)\Delta]\Big| \ .
\end{eqnarray*}
The first term on the right hand side of the above inequality
converges to zero in probability by the weak law of large numbers.
The second term
\begin{eqnarray*}
\Big|\mathbb{P}_n[ \Omega(Y) \{ \eta(Y; \theta) - \eta_0(Y;
\theta)\}\Delta ] \Big| \ \leq \  \mathbb{P}_n[\Omega(Y) \|\eta -
\eta_0\|\Delta ] \ \leq \
 K \delta_n \ \to
\ 0
\end{eqnarray*}
uniformly over $\theta$. And the last term converges uniformly to
zero in outer probability because $\{\eta_0(t; \theta): 0\leq
t\leq \tau, \theta \in \Theta\}$ is a Donsker class as we argued in
the proof of Proposition \ref{Proposition 3.1}, and both $\{\Omega(t)\}$ and $\{\Delta\}$
are also Donsker, thus $\{ \Omega(t) \eta_0(t;\theta)\Delta \}$ is Donsker and
hence a Glivenko-Cantelli class. \hspace{5mm} $\Box$

{\prop \label{Proposition 3.3} Assume the conditions in Propositions 3.1 and 3.2. Then the root of function (\ref{WCCCox}) is
asymptotically Gaussian, i.e., $n^{1/2}(\hat{\theta}_n - \theta_0)$
converges in distribution to a zero mean Gaussian random variable
with asymptotic variance $A^{-1}B\left(A^{-1}\right)'$, where
\begin{eqnarray*}
A & \!\!\!=\!\!\! & - \ \frac{\partial}{\partial \theta} \Psi(\theta,
\eta_0(\cdot;\theta))\Bigg|_{\theta=\theta_0} \, ,
\end{eqnarray*}
and
\begin{eqnarray*}
B & \!\!\!= \!\!\!& P \Bigg[ \Omega(Y) \{Z(Y) - \eta_0(Y; \theta_0)\}\Delta  \\
&& \qquad - \ \int
W(t)\{Z(t)-\eta_0(t; \theta_0)\} e^{\theta_0'Z(t)} 1(Y \geq t) d\Lambda_0(t)
\Bigg]^{\otimes 2} \, ,
\end{eqnarray*}
where $a^{\otimes 2} = aa'$.}

\smallskip

{\sc Proof:} Let ${\cal H}_0$ defined in Lemma \ref{Lemma 2.2} consist of functions of $\eta_0$ and $\hat{\eta}_n
= \hat\eta_n^{W}$, thus a Donsker class. Obviously the class of
functions $\{\psi(\theta, \eta(t;\theta)) = \Omega(t) \{Z(t) -
\eta(t;\theta)\}\Delta: \, \theta\in\Theta_0, \eta \in {\cal H}_0, 0\le t\le \tau\}$ is a
Donsker class that satisfies $P_0|\psi(\theta, \eta) - \psi(\theta_0,
\eta_0)|^2 \to 0$ as $|\theta - \theta_0| \to 0$ and $\|\eta -
\eta_0\| \to 0$ by the dominated convergence theorem. The
Fr\'{e}chet differentiability of $\{\Psi(\theta,
\eta(\cdot;\theta)): \theta \in \Theta_0, \eta \in {\cal H}_0 \}$
can be verified easily. Thus from Propositions \ref{Proposition 3.1}, \ref{Proposition 3.2} and the
remark following Lemma \ref{Lemma 2.2} together with Assumption (D), we have
all the conditions in Lemma \ref{Lemma 2.2} satisfied and thus equation
(\ref{AsymptoticRepresentation}) holds.

Now we calculate the right hand side of equation
(\ref{AsymptoticRepresentation}) for the Cox
model. Interchanging differentiation and integration yields
\begin{eqnarray*}
&& \hspace{-0.25in} n^{1/2}\dot{\Psi}_2(\theta_0, \eta_0(\cdot;\theta_0))[(\hat{\eta}_n - \eta_0)(\cdot;\theta_0)]  \\
&& \quad  = \  - \ n^{1/2}P[\Omega(Y) \{\hat{\eta}_n(Y; \theta_0) -
\eta_0(Y; \theta_0)
\} \Delta] \\
%&& \quad  = \  - \ n^{1/2}P[\{\hat{\eta}_n(Y; \theta_0) -
%\eta_0(Y; \theta_0) \} \Delta] \\
&& \quad  = \ - \ n^{1/2} \int
\Bigg[
\frac{1}{d^{(0)}(t,\theta_0)} \{D^{(1)}_n (t, \theta_0) -  d^{(1)} (t, \theta_0)\}  \\
&& \qquad \qquad    - \
\frac{D^{(1)}_n(t,\theta_0)}{D^{(0)}_n(t,\theta_0) d^{(0)}(t,
\theta_0)}\{D^{(0)}_n (t, \theta_0) - d^{(0)} (t, \theta_0)\}
\Bigg] \delta dP_{Y, \Delta}(t, \delta)  \\
&& \quad  = \  - \ n^{1/2} \int \Bigg[
\frac{1}{d^{(0)}(t,\theta_0)} \{D^{(1)}_n (t, \theta_0) -  d^{(1)} (t, \theta_0)\}  \\
&& \qquad \qquad   - \
\frac{d^{(1)}(t,\theta_0)}{d^{(0)}(t,\theta_0)^2 }\{D^{(0)}_n (t,
\theta_0) - d^{(0)} (t, \theta_0)\} \Bigg] \delta
dP_{Y, \Delta}(t, \delta) + o_{p^*}(1) \\
&& \quad  = \ - \ \mathbb{G}_n \Bigg\{\int \Big\{
W(t)\{Z(t)-\eta_0(t; \theta_0)\} e^{\theta_0'Z(t)}1(Y \geq t)\Big\} \\
&& \qquad \qquad  \times \Big\{ d^{(0)}(t, \theta_0)\Big\}^{-1}  dP_{Y,
\Delta}(t,1) \Bigg\} +  o_{p^*}(1) \, .
\end{eqnarray*}
The above second equality holds because $E(\Omega|X)=1$, and the third
equality holds because the absolute difference between the two sides except the term $o_{p^*}(1)$ becomes
\begin{eqnarray*}
&&\Bigg| \int \Bigg\{
\frac{d^{(1)}(t,\theta_0)}{d^{(0)}(t,\theta_0)^2 } -
\frac{D^{(1)}_n(t,\theta_0)}{D^{(0)}_n(t,\theta_0) d^{(0)}(t,
\theta_0)} \Bigg\} \\
&& \qquad \qquad \times \ n^{1/2}\Big \{ D^{(0)}_n(t,\theta_0) - d^{(0)}(t,
\theta_0)
\Big\}  \delta dP_{Y,\Delta}(t,\delta) \Bigg| \\
&&\qquad \leq \ \sup_{t \leq \tau} \Bigg|
\frac{d^{(1)}(t,\theta_0)}{d^{(0)}(t,\theta_0)^2 } -
\frac{D^{(1)}_n(t,\theta_0)}{D^{(0)}_n(t,\theta_0) d^{(0)}(t,
\theta_0)} \Bigg|  \\
&& \qquad \qquad \times \ \sup_{t \leq \tau} \Big| n^{1/2}\Big \{
D^{(0)}_n(t,\theta_0) - d^{(0)}(t, \theta_0) \Big\}
\Big|  \\
&&\qquad = \ o_{p^*}(1) \cdot O_{p^*}(1)    =  o_{p^*}(1)
\end{eqnarray*}
by Proposition \ref{Proposition 3.1} and tail bounds for the supremum of empirical
processes in \cite{vaart-wellner:96}, Section 2.14.

Let $G(t | \bar z(t))$ be the conditional distribution function of the
censoring time $C$ at $t$ given $\bar Z(t)= \bar z(t)$, or equivalently given $\bar Z(\tau) =\bar z(\tau)$ where $t \le \tau$, and $H_t$ be the joint distribution function
of $\bar Z(t)$. Then
\begin{eqnarray*}
d^{(0)}(t, \theta_0)  \!\!\! &=& \!\!\!  P\Big\{ W(t)1(Y \geq t)e^{\theta_0'Z(t)} \Big\}
 \\
  \!\!\! &= & \!\!\!  P\Big\{1(Y \geq t)e^{\theta_0'Z(t)} \Big\} \\
 \!\!\! &= & \!\!\! E\Big[ e^{\theta_0'Z(t)} E\{1(Y \geq t)|\bar Z(t)\}\Big] \\
 \!\!\!&  = & \!\!\!  E\Big[
e^{\theta_0'Z(t)} P(T \geq t|\bar Z(t))P(C \geq t |\bar Z(t)) \Big] \\
 \!\!\! &=&  \!\!\! \int e^{\theta_0'z(t)} {\rm exp}\Bigg\{ - \int_0^t e^{\theta_0'z(s)}d\Lambda_0(s)
\Bigg\}\{1-G(t^-|\bar z(t))\}dH_t(\bar z(t)) \, .
\end{eqnarray*}

On the other hand, from the joint distribution of $(Y, \Delta, \bar Z(Y))$, or equivalently of $(Y, \Delta, \bar Z(\tau))$,
we obtain
\begin{eqnarray*}
dP_{Y, \Delta}(t,1) &=& \Bigg[\int e^{\theta_0'z(t)} {\rm exp}\Bigg\{ - \int_0^t
e^{\theta_0'z(s)}d\Lambda_0(s) \Bigg\} \\
&& \qquad \times \{1-G(t^-|\bar z(t))\}dH_t(\bar z(t))\Bigg]
d\Lambda_0(t) \\
&=& d^{(0)}(t, \theta_0) d\Lambda_0(t).
\end{eqnarray*}
Thus we have
\begin{eqnarray*}
&& \hspace{-0.25in}  n^{1/2}\dot{\Psi}_2(\theta_0,
\eta_0(\cdot;\theta_0))[(\hat{\eta}_n - \eta_0)(\cdot;\theta_0)] \\
&&   = \,  - \, \mathbb{G}_n\Bigg[ \int \Big\{
W(t)\{Z(t)-\eta_0(t; \theta_0)\} e^{\theta_0'Z(t)}1(Y \geq t)\Big\}
d\Lambda_0(t) \Bigg] + o_{p^*}(1) .
\end{eqnarray*}
It is obvious that $\dot{\Psi}_1 = 0$, and by interchanging
differentiation and integration we have
\begin{eqnarray*}
\dot{\Psi}_2(\theta_0,
\eta_0(\cdot;\theta_0))[\dot{\eta}_0(\cdot;\theta_0)] &=& -
P\Delta\dot{\eta}_0(Y;\theta_0) \\
&=&
P\Big[\frac{\partial}{\partial\theta}\psi(\theta_0,\eta_0(\cdot;\theta))\Big]_{\theta=\theta_0}
\\
&=&
\frac{\partial}{\partial\theta}\Psi(\theta_0,\eta_0(\cdot;\theta))\Big|_{\theta=\theta_0}
\, .
\end{eqnarray*}
Then by equality
(\ref{AsymptoticRepresentation}) we have
\begin{eqnarray} \label{GeneralResultCox}
&& \hspace{-0.25in}  n^{1/2}(\hat{\theta}_n - \theta_0) \\
&&  = \ \Big\{ -
\frac{\partial}{\partial\theta}\Psi(\theta_0,\eta_0(\cdot;\theta))\Big|_{\theta=\theta_0}\Big\}^{-1}
\mathbb{G}_n \Bigg[ \Omega(Y)
\{Z(Y) - \eta_0(Y; \theta_0)\}\Delta  \nonumber\\
 && \qquad  - \ \int W(t) \{Z(t)-\eta_0(t; \theta_0)\} e^{\theta_0'Z(t)} 1(Y \geq
t) d\Lambda_0(t) \Bigg] + o_{p^*}(1) \, , \nonumber
\end{eqnarray}
which converges in distribution to a zero mean Gaussian random
variable by the central limit theorem for i.i.d. data. \hspace{5mm} $\Box$

%It can be
%easily shown that
%$$
%P\Big[\Omega\{Z - \eta_0(Y; \theta_0)\}\Delta \Big] = 0
%$$
%and
%$$ P\Bigg[\int W\{Z-\eta_0(t; \theta_0)\} e^{\theta_0'Z}1(Y
%\geq t) d\Lambda_0(t) \Bigg] = 0 \, ,
%$$
%hence the variance of $n^{1/2}(\hat{\theta}_n - \theta_0)$ can be
%estimated by $A_n^{-1}B_nA_n^{-1}$ where
%\begin{eqnarray*}
%A_n & = & - \ \frac{\partial}{\partial \theta}
%\Psi_n(\hat{\theta}_n,
%\hat{\eta}_n(\cdot;\theta))\Bigg|_{\theta=\hat{\theta}_n} \\
%B_n &=& \mathbb{P}_n \Bigg[
%\Omega \{Z - \hat{\eta}_n(Y; \hat{\theta}_n)\}\Delta \\
% && \qquad \qquad - \ \int W\{Z-\hat{\eta}_n(t; \hat{\theta}_n)\} %e^{\hat{\theta}_n'Z}1(Y \geq
%t) d\hat{\Lambda}_n(t) \Bigg]^{\otimes 2} \ ,
%\end{eqnarray*}
%where $\hat{\Lambda}_n$ can be any consistent estimator of
%$\Lambda_0$, e.g. the Breslow type of estimator similar to that in
%\cite{self-prentice:88} that only involves covariate data on
%subcohort members and failures (and hence $\hat{\Lambda}_n \to
%\Lambda_0$ weakly). \hspace{5mm} $\Box$

\medskip

It is worth noting that equation (\ref{GeneralResultCox}) reduces to the asymptotic representation of the partial likelihood estimator of \cite{cox:72}  when $\Omega_i = W_i=1$ for
all $i$. It also reduces to  the asymptotic
representation of \cite{self-prentice:88} when $\Omega_i = 1$ and $W_i$ is the inverse selection probability weight of subject $i$ into the subcohort, and of \cite{breslow-wellner:06} when $\Omega_i$ and $W_i$ are the inverse selection probability weight in a two-phase sampling design.
Note that the estimators discussed here are generally not
semiparametric efficient except the case of full-cohort data where
$\Omega_i=W_i=1$ for all $i$. Finding the most efficient estimator
is not our focus here. We refer to \cite{nan-emond-wellner:04} for
calculations of information bounds and \cite{nan:04} for an
efficient estimator when covariates are discrete.

The above calculation only considers the situation where the weights
$\Omega_i$ and $W_i$ are given for each $i$. It has been shown in the missing data literature that
using estimated rather than known weights can improve efficiency,
see e.g. \cite{rrz:94}, \cite{breslow-wellner:06}, and \cite{li-nan:11}. In particular, \cite{breslow-wellner:06} showed that, for the Cox model with time-independent covariates, the weighted estimator from a finite population sampling has the same asymptotic distribution as the weighted estimator from an i.i.d. Bernoulli sampling with the same selection probability but using the estimated weights. The asymptotic variance is smaller than that obtained using the true weights for the case of i.i.d. sampling. The same property holds for the Cox model with time-dependent covariates and time-dependent weights in the case of i.i.d. sampling. The detailed calculation follows \cite{breslow-wellner:06} and is left to the interested readers.

\subsection{Case-cohort study: the additive hazards model}

\cite{lin-ying:94} proposed the additive hazards model in which
the hazard function given covariate history $\bar Z(\cdot)$ is
$$
\lambda(t|\bar Z(t)) = \lambda_0(t) + \theta_0'Z(t) \, ,
$$
where $\lambda_0$ is the baseline hazard and $\theta_0$ is the
parameter of interest. This model allows one to estimate the covariate
effect on the absolute risk. Define the following random map:
\begin{eqnarray}
\Psi_n(\theta, \eta) &\!\!\!=\!\!\!& \frac{1}{n}\sum_{i=1}^n \Bigg\{\Omega_i(Y_i) \{Z_i(Y_i) -
\eta(Y_i)\}\Delta_i  \label{AdditiveHazards} \\
&& \qquad \quad - \ \int \Omega_i(t) \{Z_i(t) - \eta(t)\} 1(Y_i \geq t)
\theta' Z_i(t) \, dt \Bigg\} \nonumber
\end{eqnarray}
with
$$
\eta_0(t) = \frac{E\{Z(t)1(Y \geq t)\}}{E\{1(Y \geq t)\}} \, ,
$$
where $\Omega_i$ are defined in the same way
as that in the previous subsection for the Cox model. Then the
estimating function proposed by \cite{lin-ying:94} can be viewed
as the above function (\ref{AdditiveHazards}) with $\Omega_i = 1$
and $\eta_0$ being estimated empirically, which has the following
form:
\begin{eqnarray}
\Psi_n(\theta, \tilde\eta_n^F) &\!\!\!=\!\!\!& \frac{1}{n}\sum_{i=1}^n \Bigg\{ \{Z_i(Y_i) -
\tilde\eta_n^F(Y_i)\}\Delta_i \label{FullAdditiveHazards} \\
&& \qquad \quad - \ \int \{Z_i(t) - \tilde\eta_n^F(t)\} 1(Y_i \geq t)
\theta' Z_i(t) \, dt \Bigg\}  \nonumber
\end{eqnarray}
with
$$
\tilde\eta_n^F(t) = \frac{\sum_{j=1}^n Z_j(t) 1(Y_j \geq t)}{\sum_{j=1}^n 1(Y_j
\geq t)} \ .
$$
Note that both $\eta_0$ and $\tilde\eta_n^F$ do not involve $\theta$.
The estimator of $\theta$ has an explicit form:
\begin{eqnarray}\label{ThetaHatAdditive}
\tilde{\theta}_n = \Bigg[\frac{1}{n}\sum_{i=1}^n \int \{Z_i(t) -
\tilde\eta_n^F(t) \}^{\otimes 2} 1(Y_i \geq t) dt \Bigg]^{-1}  \frac{1}{n}\sum_{i=1}^n  \{Z_i(Y_i) - \tilde\eta_n^F(Y_i)\}\Delta_i  .
\end{eqnarray}
\cite{lin-ying:94} defined the above $\Psi_n(\theta, \tilde\eta_n^F)$ and
$\tilde{\theta}_n$ using the stochastic integral formulation and
studied their asymptotic properties using martingale theory.

For case-cohort studies, \cite{kulich-lin:00} modified the estimating
function (\ref{FullAdditiveHazards}) and proposed the following
estimating function (with $\Omega_i = W_i$):
\begin{eqnarray}
\qquad \Psi_{n}(\theta, \tilde\eta_n^{W}) &\!\!\!=\!\!\!& \frac{1}{n}\sum_{i=1}^n
\Bigg\{\Omega_i(Y_i) \{Z_i(Y_i) - \tilde\eta_n^{W}(Y_i)\}\Delta_i \label{WCCAdditiveHazards}
\\
&& \qquad  - \ \int \Omega_i(t) \{Z_i(t) - \tilde\eta_n^{W}(t)\} 1(Y_i \geq t)
\theta' Z_i(t) \, dt \Bigg\} \nonumber
\end{eqnarray}
with
\begin{eqnarray}
&& \tilde\eta_n^{W}(t) = \frac{\sum_{j=1}^n W_j(t) Z_j(t) 1(Y_j \geq t)}{\sum_{j=1}^n W_j(t) 1(Y_j \geq
t)} \ .
\label{WCCAdditiveHazardsEta}
\end{eqnarray}

The estimator again has an explicit form
\begin{eqnarray}
\tilde{\theta}_n &\!\!\!=\!\!\!& \Bigg[\frac{1}{n}\sum_{i=1}^n \int \Omega_i(t)
\{Z_i(t) - \tilde\eta_n^{W}(t) \} Z_i(t)' 1(Y_i \geq t) dt \Bigg]^{-1} \label{ThetaTildeAdditive} \\
&& \qquad \qquad \times
\frac{1}{n}\sum_{i=1}^n \Omega_i(Y_i) \{Z_i(Y_i) -
\tilde\eta_n^{W}(Y_i)\}\Delta_i \ .  \nonumber
\end{eqnarray}
Here we have extended the method of \cite{kulich-lin:00} by
introducing two weight matrices $\Omega$ and $W$ in
(\ref{AdditiveHazards}) and (\ref{WCCAdditiveHazardsEta}),
respectively, as in the previous subsection.

When weights $W_i$ or $\Omega_i$ depend on $\Delta_i$ as in
(\ref{Weights}), for the same reason as that in the previous
example,
martingale theory does not apply. Here we provide a
proof without using stochastic
integrals. As we assumed for the Cox model, $\Omega_i$ and $W_i$ are
nonnegative with unit conditional expectation given
complete data $X_i$.

We consider the weighted estimating function
(\ref{WCCAdditiveHazards}) that reduces to
(\ref{FullAdditiveHazards}) when $\Omega_i = W_i =1$ for all $i$. Without loss of generality, we
assume one-dimensional covariate $Z$ and thus one-dimensional
$\theta$ in the following calculation. Multi-dimensional case is a
straightforward extension.

{\prop \label{Proposition 3.4} Let $\tilde{\eta}_n(t) =
\hat\eta_n^{W}(t)$ as in equation (\ref{WCCAdditiveHazards}). Suppose the weight process $W(t)$ has bounded sample paths of bounded variation. Then both $\tilde\eta_n(t)$ and $\eta_0(t)$ belong to a Donsker class, and further  we
have $\|\tilde{\eta}_n - \eta_0 \| = O_{p^*}(n^{-1/2})$. }

%\smallskip

{\sc Proof:} 
%The function $\tilde\eta_n^{W}(t)$  is exactly the same as that in equation
%(\ref{WCCCox}) at $\theta = 0$, and $\eta_0(t)$ is the same as its
%counterpart in the Cox model at $\theta=0$. Hence the stated result
This is a direct consequence of Proposition \ref{Proposition 3.1} with $\theta=0$. \hspace{5mm} $\Box$

{\prop \label{Proposition 3.5} Assume the conditions in Proposition 3.4 and suppose the weight process $\Omega(t)$ also has bounded sample paths of bounded variation. Then the root of function
(\ref{WCCAdditiveHazards}) is a consistent estimator of $\theta_0$.
}

\smallskip

{\sc Proof:} Similar to the proof of Proposition \ref{Proposition 3.2}, we only need
to verify those conditions in Lemma \ref{Lemma 2.1}. Obviously $\Psi(\theta,
\eta_0) = P\{\psi(\theta,\eta_0)\}$ is a linear function for
$\theta$ with a non-zero slope by Assumption (D), hence $\theta_0$ is the unique solution of $\Psi(\theta,
\eta_0) = 0$. Proposition \ref{Proposition 3.4} provides the uniform
consistency of $\hat{\eta}_n$.
We now verify condition (\ref{ConsistencyCondition}). Let $\| \eta -
\eta_0 \| \downarrow 0$. We have
\begin{eqnarray*}
&& \hspace{-0.25in} |\Psi_n(\theta, \eta) - \Psi(\theta, \eta_0)|  \\
&&  \le \ \big| \mathbb{P}_n[\Omega(Y) \{ Z(Y) - \eta(Y)\}\Delta] -
P[\Omega(Y)\{Z(Y)-\eta_0(Y)\}\Delta]\big| \\
&&  \qquad + \ \bigg|\mathbb{P}_n \int \Omega(t) \{Z(t) - \eta(t)\}1(Y\geq
t)\theta Z(t) \ dt \\
&&   \qquad \qquad - \ P \int \Omega(t)\{Z(t)-\eta_0(t)\}1(Y\geq t)\theta Z(t) \ dt
\bigg| \\
&&  \le \ \big| (\mathbb{P}_n-P)[\Omega(Y)\{Z(Y)-\eta_0(Y)\}\Delta]\big| + \big|\mathbb{P}_n[\Omega(Y)\{\eta(Y)-\eta_0(Y)\}\Delta]\big| \\
%\end{eqnarray*}
%The uniform convergence to zero in probability for the first term on
%the right hand side of the above inequality can be shown by similar
%arguments as in the proof of Proposition \ref{Proposition 3.2}. The second %term can be
%written as
%\begin{eqnarray*}
&& \qquad + \ \bigg|(\mathbb{P}_n - P) \int
\Omega(t)\{Z(t)-\eta_0(t)\}1(Y\geq t)\theta Z(t) \, dt \bigg|\\
&& \qquad \qquad + \
\bigg|\mathbb{P}_n \int \Omega(t)\{\eta(t)-\eta_0(t)\}1(Y\geq t)\theta Z(t) \ dt \bigg| \\
&&  \leq \ \big| (\mathbb{P}_n-P)[\Omega(Y)\{Z(Y)-\eta_0(Y)\}\Delta]\big| \\
&& \qquad + \ \bigg|(\mathbb{P}_n - P) \int
\Omega(t)\{Z(t)-\eta_0(t)\}1(Y\geq t)\theta Z(t) \, dt \bigg| \\
&& \qquad \qquad  + \ \delta_n
\mathbb{P}_n \bigg\{\Omega(Y)\Delta + \int \Omega(t) 1(Y \geq t) |\theta Z(t)| \, dt \bigg\},
\end{eqnarray*}
in which the first two terms on the right hand side of the last inequality  converge to zero in
probability by the weak law of large
numbers, and the third term converges to zero because $\delta_n  \to 0$. We then have
the desired result by Lemma \ref{Lemma 2.1}.
%Note that the result also holds
%when $|\theta| = O(n^{1/2-\epsilon})$ for any $\epsilon
%> 0$ by Proposition \ref{Proposition 3.4}.
\hspace{5mm} $\Box$

{\prop \label{Proposition 3.6} Assume the conditions in Propositions 3.4 and 3.5. Then the root of function
(\ref{WCCAdditiveHazards}), given in (\ref{ThetaTildeAdditive}), is asymptotically Gaussian, i.e.,
$n^{1/2}(\tilde{\theta}_n - \theta_0)$ converges in distribution
to a zero mean Gaussian random variable. }

\smallskip

{\sc Proof:} The proof can proceed either from (\ref{ThetaTildeAdditive}) directly or by using Corollary \ref{Corollary 2.3}. We show the latter. Similar to the proof
of Proposition 3.3, the Fr\'{e}chet differentiability of
$\{\Psi(\theta, \eta): \theta \in \Theta_0, \eta \in {\cal H}_0 \}$
can be verified easily. Obviously the set $\{ \Omega(t) \Delta\{Z(t) - \eta(t)\}:  \eta \in {\cal H}_0, 0\le t \le \tau \}$ is Donsker, thus we only need to show the class of
functions $\{ \int_0^\tau \Omega(t) \{Z(t) - \eta(t)\}1(Y \geq t) \theta Z(t) dt:
\theta \in \Theta_0, \eta \in {\cal H}_0\}$ is Donsker, here ${\cal
H}_0$ is reduced from that in the proof of Proposition 3.3.
% Now
%\begin{eqnarray*}
%&& \int_0^\tau \Omega(t) \{Z(t) - \eta(t)\}1(Y \geq t) \theta Z(t) dt \\
%&& \qquad = \theta
%\int_0^\tau\Omega(t) Z(t)^2 1(Y \geq t) dt - \theta \int_0^\tau \Omega(t) %Z(t) \eta(t)1(Y \geq t)dt \, .
%\end{eqnarray*}
%and both $\Omega(t)$ and $\Omega(t)\eta(t)$ belong to Donsker
%classes, we only need to show $f_\eta(Y) \equiv \int\eta(t) 1(Y \geq
%t)dt$ forms a Donsker class.
Let $f = \int_0^\tau \Omega(t) \{Z(t) - \eta(t)\} Z(t)1(Y \geq t) dt$ and
$$
f^m = \sum_{i=1}^m \Omega(t_i) \{Z(t_i) - \eta(t_i)\}Z(t_i) 1(Y \geq t_i) (t_{i+1}-t_i) =
\sum_{i=1}^m f_i \lambda_i \, ,
$$
where
$$
f_i = \Omega(t_i) \{Z(t_i) - \eta(t_i)\}Z(t_i) 1(Y \geq t_i) , \quad \lambda_i = t_{i+1}-t_i,
$$
and $\{(t_1, t_2], \dots, (t_m, \tau]\}$ forms a partition of the
interval $(0,\tau]$. The set $\{ f^m\}$ is the convex hull of ${\cal
F} = \{f_i\}$, and thus a Donsker class by Theorem 2.10.3 in
\cite{vaart-wellner:96} since ${\cal F} $ is Donsker. Now we know
that $f^m \to f$ both pointwise and in $L_2(P)$ by the
boundedness of $Y$ and $\eta$, then $\{f(\cdot) \}$ is Donsker by
Theorem 2.10.2 in \cite{vaart-wellner:96}.

We then calculate the right hand side of equation
(\ref{AsymptoticRepresentation2}). Direct calculation yields
\begin{eqnarray}
\hspace{-0.25in}  n^{1/2} \dot{\Psi}_\eta(\theta_0,
\eta_0)(\tilde{\eta}_n -
\eta_0)  &\!\!\!=\!\!\!&  - \ n^{1/2}P[\{\tilde{\eta}_n(Y) - \eta_0(Y)\}\Delta] \label{psi-eta}\\
&&    + \
n^{1/2}P\bigg[\int \{\tilde{\eta}_n(t) - \eta_0(t)\}1(Y \geq t)
\theta _0 Z(t) dt \bigg] \, \nonumber 
\end{eqnarray}
by applying $E(\Omega|X) = 1$. Let $d^{(0)}(t) \equiv P\{W(t)1(Y \geq t)\}
= P\{1(Y \geq t)\}$ and $d^{(1)}(t) \equiv P\{W(t)Z(t)1(Y \geq t)\} = P\{Z(t)1(Y
\geq t)\}$, where $E(W|X) = 1$. Similar to the proof of Proposition 3.3, the first term on the right hand side
of equation (\ref{psi-eta}) can be written as
\begin{eqnarray*}
&& - \ n^{1/2}P[\{\tilde{\eta}_n(Y) - \eta_0(Y)\}\Delta] \\
&& \qquad  = \ - \ \mathbb{G}_n\bigg[ \int W(t)
\{Z(t) - \eta_0(t)\} 1(Y \geq t)  d^{(0)}(t)^{-1} dP_{Y, \Delta}(t,1)\bigg]  + o_{p^*}(1)
\\
&& \qquad  = \ - \ \mathbb{G}_n \bigg[ \int W(t) \{Z(t) -
\eta_0(t) \} 1(Y \geq t) \lambda_0(t) dt  \\
&& \qquad \qquad    + \  \int W(t) \{Z(t) -
\eta_0(t)\} 1(Y \geq t)\theta_0  \eta_0(t) dt \bigg] + o_{p^*}(1)
\end{eqnarray*}
since from the joint distribution of $(Y, \Delta, \bar Z(Y))$ we have
\begin{eqnarray*}
\frac{dP_{Y,\Delta}(t,1)}{dt} &=& \int \{\lambda_0(t) + \theta_0 z(t)\}
\{ 1 -
F(t|\bar z(t))\}  \\
&& \quad \{1- G(t^-|\bar z(t))\} dH_t(\bar z(t)) \\
&=& \lambda_0(t) P\{1(Y \geq t)\} + \theta_0 P\{Z(t)1(Y \geq t)\} \\
&=& \lambda_0 (t) d^{(0)}(t) + \theta_0 d^{(1)}(t) \, .
\end{eqnarray*}
From the proof of Proposition \ref{Proposition 3.1} we have
\begin{eqnarray*}
n^{1/2}\{ \tilde{\eta}_n(t) - \eta_0(t) \} = d^{(0)}(t)^{-1}
\mathbb{G}_n\Big[W(t)\{Z(t) - \eta_0(t)\}1(Y \geq t)\Big] + o_{p^*}(1) \,
,
\end{eqnarray*}
so the second term on the right hand side of (\ref{psi-eta}) can be
rewritten as
\begin{eqnarray*}
&& \!\!\!\!\!\! \int n^{1/2} \{\tilde{\eta}_n(t) - \eta_0(t)
\} P\{1(Y \geq t) \theta_0 Z(t)\}  \, dt \nonumber \\
&& \quad  = \ \int d^{(0)}(t)^{-1}  \mathbb{G}_n \Big[W(t)\{Z(t)
-\eta_0(t)\}1(Y \geq
t)\Big]  \theta_0 d^{(1)}(t) \, dt + o_{p^*}(1) \nonumber \\
&& \quad = \ \mathbb{G}_n \bigg [ \int W(t) \{Z(t) -
\eta_0(t)\} 1(Y \geq t) \theta_0  \eta_0(t) \, dt\bigg] + o_{p^*}(1) \, .
\end{eqnarray*}
Thus from (\ref{AsymptoticRepresentation2}) we obtain
\begin{eqnarray} \label{GeneralResultAdditive}
n^{1/2}(\tilde{\theta}_n - \theta_0) &\!\!\!=
\!\!\! & \bigg[P \bigg\{ \int \Omega(t) \{Z(t) - \eta_0(t)\}1(Y \geq t) Z(t) \, dt
\bigg\}\bigg]^{-1}  \\
&& \quad \times \ \mathbb{G}_n \bigg[\Omega(Y) \{Z(Y)-\eta_0(Y)\}\Delta  \nonumber \\
&& \quad - \ \int \{ \Omega(t)\theta_0 Z(t) + W(t) \lambda_0(t)\} \{Z(t)-\eta_0(t)\}
1(Y\geq t) \, dt \bigg] \nonumber \\
&& \quad + \ o_{p^*}(1)  , \nonumber
\end{eqnarray}
which is asymptotic normal by the central limit theorem. This
asymptotic representation reduces to that in \cite{kulich-lin:00}
when $\Omega_i = W_i$.  Again, we
do not require $\Omega_i$ and $W_i$ to be predictable.

%The variance estimator of
%$\tilde{\theta}_n$ can be obtained from the above expression by
%replacing $\theta_0$, $\eta_0$, and the baseline cumulative hazard
%$\Lambda_0(t)$ by their estimators. We refer to
%\cite{kulich-lin:00} for an estimator of $\Lambda_0(t)$.

\section{Discussion}

We consider i.i.d. sampling for the case-cohort studies. \cite{breslow-wellner:06} have considered finite population stratified sampling and applied the exchangeably weighted bootstrap empirical process theory of \cite{p-wellner:93} for the Cox model with time-independent covariates. The general Z-estimation theory in Section 2 is likely to be applicable to the finite population stratified sampling designs for time-dependent covariates.

The theory in Section 2 requires smooth $\eta$ with respect to $\theta$, which is mainly restricted by the smoothness condition (\ref{hu:smooth}) or (\ref{huang:smooth}). For non-smooth $\eta$, for example, the rank-based estimating function for the accelerated failure time model, the smoothness condition does not hold. \cite{nan-kalbfleisch-yu:09} have showed that a similar idea for bundled parameters with missing data is applicable to the rank-based estimator for the accelerated failure time model. For models with bundled parameters in the original parameterization, \cite{ding-nan:11} have proposed a sieve maximum likelihood estimating method and applied the method to the efficient estimation of the accelerated failure time model.

We have discussed two examples, the proportional hazards model and
the additive hazards model in case-cohort studies, though our method
applies to a much broader range of semiparametric estimation
problems. The parameter estimation in the case-cohort studies is
hard to handle by traditional martingale based methods when certain
more efficient but unpredictable weights are considered, but becomes
straightforward by using the general pseudo Z-estimation theory.

Another point worth mentioning is that for missing data problems, the
estimated likelihood method of \cite{pepe-fleming:91}, the mean
score method of \cite{reilly-pepe:95}, and the pseudoscore method
of \cite{chatterjee-chen-breslow:03}, among others, also fit into
the general Z-estimation framework nicely. Let $Y$ be the response
variable and $(Z, V)$ be covariates where $Z$ can be missing
sometimes. Let $R$ be the indicator that takes value 1 if $Z$ is
observed and $0$ otherwise. Let $X$ denote the observed data.
Suppose that the parameter of interest $\theta \in \Theta \subset
\mathbb{R}^d$ could be estimated by using the complete data score
function $\dot{l}_\theta^0(\cdot; \theta)$ as the estimating
function if there were no missing data. When $Z$ is sometimes
missing at random \citep{little-rubin:02}, then the observed
data score function for $\theta$ becomes
\begin{eqnarray*}
\dot{l}_\theta(X; \theta, \eta_0(\cdot; \theta)) =  R
\dot{l}^0_\theta(Y,Z,V; \theta) + (1-R) \eta_0(Y,V; \theta) \, ,
\end{eqnarray*}
where $\eta_0(Y,V; \theta) = E\{ \dot{l}^0_\theta(Y,Z,V;\theta) | Y,
V\}$ whose functional form is unknown. Define $\psi(\cdot; \theta,
\eta(\cdot;\theta)) = \dot{l}_\theta(\cdot; \theta,
\eta(\cdot;\theta))$. Then $\psi(\cdot; \theta,
\hat{\eta}_n(\cdot;\theta))$ becomes an estimating function
for $\theta$ where $\hat{\eta}_n(\cdot;\theta))$ is an estimator of
$\eta_0(\cdot;\theta)$. The asymptotic properties of the
Z-estimator for $\theta$ depend on the behavior of $\hat{\eta}_n$
and may be derived from the theorems given in Section 2. Authors of
aforementioned references have proposed nonparametric methods to
estimate $\eta_0(\cdot; \theta)$.
Apparently efficiency can be improved by using the weighted
estimating function proposed by \cite{rrz:94}. The proposed methodology may also apply to the composite likelihoods for semiparametric models, see e.g. \cite{lindsay:87} and \cite{ varin-etal:11}, particularly for missing data problems.

\section*{Appendix: Proofs of Lemmas \ref{Lemma 2.1} and \ref{Lemma 2.2}}

\subsection*{Proof of Lemma \ref{Lemma 2.1}} Since $\theta_{0}$ is the unique
solution to $\Psi(\theta, \eta_{0}(\cdot;\theta))=0$, this implies
that for any fixed $\epsilon>0$, there exists a $\delta>0$ such that
$$P\left[|\hat{\theta}_{n}-\theta_{0}|>\epsilon\right]
\leq P\left[|\Psi(\hat{\theta}_{n}, \eta_{0}(\cdot;
\hat{\theta}_n))|> \delta\right].$$ If we can prove
$|\Psi(\hat{\theta}_{n}, \eta_{0}(\cdot;
\hat{\theta}_n))|\rightarrow_{p^{*}} 0$, then the consistency of
$\hat{\theta}_{n}$ will follow immediately.

To do this, first note that since
$||\hat{\eta}_{n}-\eta_{0}||=o_{p^{*}}(1)$, there exists a sequence
$\{\delta_{n}\}\downarrow 0 $ such that
$||\hat{\eta}_{n}-\eta_{0}||\leq \delta_{n}$ with probability
tending to one. Hence taking $\eta=\hat{\eta}_{n}$ in equation
(\ref{ConsistencyCondition}), we have the following inequalities:
\begin{eqnarray*}
|\Psi(\hat{\theta}_{n}, \eta_{0}(\cdot;\hat{\theta}_n))| \!\!\!&\leq&
\!\!\! |\Psi_{n}(\hat{\theta}_{n},
\hat{\eta}_{n}(\cdot;\hat{\theta}_n))| + |\Psi(\hat{\theta}_{n},
\eta_{0}(\cdot;\hat{\theta}_n))
-\Psi_{n}(\hat{\theta}_{n}, \hat{\eta}_{n}(\cdot;\hat{\theta}_n))|\\
\!\!\!&\leq& \!\!\! |\Psi_{n}(\hat{\theta}_{n},
\hat{\eta}_{n}(\cdot;\hat{\theta}_n))| + o_{p^{*}}\Big( 1 +
|\Psi_{n}(\hat{\theta}_{n},\hat{\eta}_{n}(\cdot;\hat{\theta}_n))| \\
&& \qquad + \
|\Psi(\hat{\theta}_{n}, \eta_{0}(\cdot;\hat{\theta}_n))|\Big) \\
\!\!\!&\leq& \!\!\! o_{p^{*}}(1) + o_{p^{*}}\left(1+
o_{p^{*}}\left(1\right) + |\Psi(\hat{\theta}_{n},
\eta_{0}(\cdot;\hat{\theta}_n))|\right),
\end{eqnarray*}
which implies $|\Psi(\hat{\theta}_{n},
\eta_{0}(\cdot;\hat{\theta}_n))|=o_{p^{*}}(1)$. So we have proved
the consistency of pseudo Z-estimators $\hat{\theta}_{n}$.
\hspace{0.5cm}$\Box$

\subsection*{Proof of Lemma \ref{Lemma 2.2}} We first show a result that we will use
in the proof: under Conditions (i) and (ii),
\begin{eqnarray} \label{import}
n^{1/2}\left|\Psi(\hat{\theta}_{n},\hat{\eta}_{n}(\cdot,
\hat{\theta}_n))\right|=O_{p^{*}}(1).
\end{eqnarray}
By Condition (i), we have the following inequality:
\begin{eqnarray*}
&&\!\!\!\!\!\! n^{1/2}\left|(\Psi_{n}-\Psi)(\hat{\theta}_{n},
\hat{\eta}_{n}(\cdot;\hat{\theta}_n))-
(\Psi_{n}-\Psi)(\theta_{0},\eta_{0}(\cdot;\theta_0))\right|\\
  &&  = \
o_{p^{*}}(1)+o_{p^{*}}\left(n^{1/2}\left|\Psi_{n}(\hat{\theta}_{n},\hat{\eta}_{n}(\cdot;\hat{\theta}_n))\right|\right) +
o_{p^{*}}\left(n^{1/2}\left|\Psi(\hat{\theta}_{n},\hat{\eta}_{n}(\cdot;\hat{\theta}_n))\right|\right).
\end{eqnarray*}
By the triangle inequality $-|a|+|b|-|c|\leq |a-b-c|$ and the fact
that $\Psi(\theta_{0},\eta_{0}(\cdot;\theta_0))=0$,
\begin{eqnarray*}
&&\hspace{-2cm}n^{1/2}\left|\Psi(\hat{\theta}_{n},\hat{\eta}_{n}(\cdot;\hat{\theta}_n))
\right|-n^{1/2}\left|\Psi_{n}(\hat{\theta}_{n},
\hat{\eta}_{n}(\cdot;\hat{\theta}_n))\right|-n^{1/2}\left|\Psi_{n}(\theta_{0},\eta_{0}(\cdot;\theta_0))\right|\\
\qquad \qquad &\leq&n^{1/2}\left|(\Psi_{n}-\Psi)(\hat{\theta}_{n},
\hat{\eta}_{n}(\cdot;\hat{\theta}_n))-(\Psi_{n}-\Psi)
(\theta_{0},\eta_{0}(\cdot;\theta_0))\right|\\
\qquad \qquad
&=&o_{p^{*}}(1)+o_{p^{*}}\left(n^{1/2}\left|\Psi_{n}(\hat{\theta}_{n},\hat{\eta}_{n}(\cdot;\hat{\theta}_n))\right|\right)
\\
&& \qquad + \ o_{p^{*}}\left(n^{1/2}\left|\Psi(\hat{\theta}_{n},\hat{\eta}_{n}(\cdot;\hat{\theta}_n))\right|\right),
\end{eqnarray*}
which implies
\begin{eqnarray*}
&&\hspace{-2cm}n^{1/2}\left|\Psi(\hat{\theta}_{n},\hat{\eta}_{n}(\cdot;\hat{\theta}_n))\right|[1-o_{p^{*}}(1)]\\
\qquad \qquad &\leq&
o_{p^{*}}(1)+n^{1/2}\left|\Psi_{n}(\hat{\theta}_{n},\hat{\eta}_{n}(\cdot;\hat{\theta}_n))\right|[1+o_{p^{*}}(1)]
\\
&& \qquad + n^{1/2}\left|\Psi_{n}(\theta_{0},\eta_{0}(\cdot;\theta_0))\right|\\
\qquad \qquad &=& o_{p^{*}}(1) + o_{p^{*}}(1) + O_{p^{*}}(1).
\end{eqnarray*}
Hence (\ref{import}) holds.

We then show the root-$n$ consistency of $\hat{\theta}_n$. Since
$|\hat{\theta}_{n}-\theta_{0}|=o_{p^{*}}(1)$ and
$||\hat{\eta}_{n}-\eta_{0}||=O_{p^{*}}(n^{-\beta})$ with $\beta>0$,
there exists a sequence $\{\delta_{n}\}\downarrow 0 $ and $c>0$ such
that $|\hat{\theta}_{n}-\theta_{0}|\leq\delta_{n}$ and
$||\hat{\eta}_{n}-\eta_{0}||\leq cn^{-\beta}$ with probability
approaching one. Hence taking $(\theta,
\eta)=(\hat{\theta}_{n},\hat{\eta}_{n})$ in the smoothness condition
(\ref{huang:smooth}):
\begin{eqnarray}
\label{pf1} \nonumber && \hspace{-1cm}
\left|n^{1/2}\left\{\Psi(\hat{\theta}_{n},
\hat{\eta}_{n}(\cdot;\hat{\theta}_n))- \Psi(\theta_{0},
\eta_{0}(\cdot;\theta_0))\right\}\right.
\\
&& \qquad - \
n^{1/2}\left\{\dot{\Psi}_1(\theta_0,\eta_0(\cdot;\theta_0))+\dot{\Psi}_2(\theta_0,
\eta_0(\cdot;\theta_0))[\dot{\eta}_0(\cdot;\theta_0)]\right\}(\hat{\theta}_{n}-\theta_{0})
  \nonumber \\
&& \qquad  \left. - \ n^{1/2}\dot{\Psi}_2(\theta_{0},
\eta_{0}(\cdot;\theta_0))
[(\hat{\eta}_{n}-\eta_{0})(\cdot;\theta_0)] \right| \nonumber \\
&&  = \ o_{p^{*}}\left(n^{1/2}|\hat{\theta}_{n}-\theta_{0}|\right)+
O_{p^{*}}\left(n^{1/2}\|\hat{\eta}_{n}-\eta_{0}\|^{\alpha} \right)\nonumber \\
&&  = \ o_{p^{*}}\left(1+n^{1/2}|\hat{\theta}_{n}-\theta_{0}|\right),
\end{eqnarray}
since
$n^{1/2}O_{p^{*}}(||\hat{\eta}_{n}-\eta_{0}||^{\alpha})=o_{p^{*}}(1)$
by $\alpha\beta>1/2$. Same result can be obtained by using the
smoothness condition (\ref{hu:smooth}) for $\beta=1/2$. By equation
(\ref{import}), the fact that
$\Psi(\theta_{0},\eta_{0}(\cdot;\theta_0))=0$, and the triangle
inequality $-|a|+|b|-|c|\leq |a-b-c|$, equation (\ref{pf1}) implies
\begin{eqnarray}
\nonumber &&
\hspace{-1.0cm}-O_{p^{*}}(1)+\left|n^{1/2}\left\{\dot{\Psi}_1(\theta_0,\eta_0(\cdot;\theta_0))+\dot{\Psi}_2(\theta_0,
\eta_0(\cdot;\theta_0))[\dot{\eta}_0(\cdot;\theta_0)]\right\}(\hat{\theta}_{n}-\theta_{0})\right| \\
&& \qquad   - \ \left|n^{1/2}\dot{\Psi}_2(\theta_{0},
\eta_{0}(\cdot;\theta_0))[(\hat{\eta}_{n}
-\eta_{0})(\cdot;\theta_0)] \right| \nonumber \\
&& \leq \
o_{p^{*}}\left(1+n^{1/2}\left|\hat{\theta}_{n}-\theta_{0}\right|\right)
. \label{rateproof}
\end{eqnarray}
Since the $d\times d$ matrix
$\dot{\Psi}_1(\theta_0,\eta_0(\cdot;\theta_0))+\dot{\Psi}_2(\theta_0,
\eta_0(\cdot;\theta_0))[\dot{\eta}_0(\cdot;\theta_0)]$ is
nonsingular, there exist a constant $c_{1}>0$ such that
$$\left|\left\{\dot{\Psi}_1(\theta_0,\eta_0(\cdot;\theta_0))+\dot{\Psi}_2(\theta_0,
\eta_0(\cdot;\theta_0))[\dot{\eta}_0(\cdot;\theta_0)]\right\}(\theta-\theta_{0})\right|\geq
c_{1}|\theta-\theta_{0}|$$ for $|\theta-\theta_{0}|\rightarrow 0$. On
the other hand, by Condition (iv), combination with inequality
(\ref{rateproof}) yields
\begin{eqnarray*}
O_{p^{*}}(1)&\geq&
\left|n^{1/2}\left\{\dot{\Psi}_1(\theta_0,\eta_0(\cdot;\theta_0))+\dot{\Psi}_2(\theta_0,
\eta_0(\cdot;\theta_0))[\dot{\eta}_0(\cdot;\theta_0)]\right\}(\hat{\theta}_{n}-\theta_{0})\right|
\\
&& \qquad - \, \left|n^{1/2}\dot{\Psi}_2(\theta_{0},
\eta_{0}(\cdot;\theta_0))[(\hat{\eta}_{n}
-\eta_{0})(\cdot;\theta_0)] \right|  \\
&& \qquad - \ o_{p^{*}}\left(1+n^{1/2}\left|\hat{\theta}_{n}-\theta_{0}\right|\right)\\
&\geq&
c_{1}n^{1/2}\left|\hat{\theta}_{n}-\theta_{0}\right|-O_{p^*}(1)-o_{p^{*}}\left(1+n^{1/2}\left|\hat{\theta}_{n}
-\theta_{0}\right|\right)\\
&=&
\left\{O_{p^{*}}(1)-o_{p^{*}}(1)\right\}n^{1/2}\left|\hat{\theta}_{n}-\theta_{0}\right|-O_{p^{*}}(1).
\end{eqnarray*}
Hence the sequence $n^{1/2}\left|\hat{\theta}_{n}-\theta_{0}\right|$
must be bounded in outer probability.

Now we are ready to prove equation (\ref{AsymptoticRepresentation}).
Because
\begin{eqnarray}
\nonumber
&&\hspace{-0.5cm}n^{1/2}\left[\Psi(\hat{\theta}_{n}, \hat{\eta}_{n}(\cdot;\hat{\theta}_n))-\Psi(\theta_{0}, \eta_{0}(\cdot;\theta_0))\right]\\
\nonumber && \qquad = \ n^{1/2}\left[\Psi(\hat{\theta}_{n},
\hat{\eta}_{n}(\cdot;\hat{\theta}_n))-\Psi_{n}(\hat{\theta}_{n},
\hat{\eta}_{n}(\cdot;\hat{\theta}_n)) \right. \nonumber \\
&& \qquad \qquad + \left. \Psi_{n}(\hat{\theta}_{n}, \hat{\eta}_{n}(\cdot;\hat{\theta}_n))-\Psi(\theta_{0}, \eta_{0}(\cdot;\theta_0))\right] \nonumber \\
\nonumber
&& \qquad = \ n^{1/2}(\Psi-\Psi_{n})(\hat{\theta}_{n},\hat{\eta}_{n}(\cdot;\hat{\theta}_n))+o_{p^{*}}({
1})- 0 \\
\nonumber && \qquad = \ - \ n^{1/2}(\Psi_{n}-\Psi)(\theta_{0},
\eta_{0}(\cdot;\theta_0))\pm
o_{p^{*}}\Big(1+n^{1/2}\left|\Psi_{n}(\hat{\theta}_{n},
\hat{\eta}_{n}(\cdot;\hat{\theta}_n))\right| \nonumber \\
&& \qquad \qquad \qquad  + \ n^{1/2}\left|\Psi(\hat{\theta}_{n},\hat{\eta}_{n}(\cdot;\hat{\theta}_n))\right|\Big) \hspace{5mm}(\mbox{by Condition (i)}) \nonumber\\
&& \qquad = \ - \ n^{1/2}(\Psi_{n}-\Psi)(\theta_{0}, \eta_{0}(\cdot;\theta_0))
\pm o_{p^{*}}(1)\hspace{2mm}(\mbox{by equation (\ref{import})}),
\label{pf2}
\end{eqnarray}
after replacing equation (\ref{pf2}) into the first term in the
first line of equation (\ref{pf1}) we obtain
\begin{eqnarray*}
&&\hspace{-1cm}\left|-n^{1/2}(\Psi_{n}-\Psi)(\theta_{0},
\eta_{0}(\cdot;\theta_0))\pm
o_{p^{*}}(1)-n^{1/2}\left\{\dot{\Psi}_1(\theta_0,\eta_0(\cdot;\theta_0))
\right. \right. \\
&& \qquad  \left.+ \ \dot{\Psi}_2(\theta_0,
\eta_0(\cdot;\theta_0))[\dot{\eta}_0(\cdot;\theta_0)]\right\}(\hat{\theta}_{n}-\theta_{0}) \\
&& \qquad  \left. - \
n^{1/2}\dot{\Psi}_2(\theta_{0}, \eta_{0}(\cdot;\theta_0))
[(\hat{\eta}_{n}-\eta_{0})(\cdot;\theta_0)] \right| \nonumber \\
 &&  = \ o_{p^{*}}\left(1+n^{1/2}\left|\hat{\theta}_{n}-\theta_{0}\right|\right)\\
\nonumber  &&  = \ o_{p^{*}}(1),
\end{eqnarray*}
which implies
\begin{eqnarray*}
n^{1/2}(\hat{\theta}_n - \theta_0) &\!\!\!=\!\!\!& \Big\{ -
\dot{\Psi}_1(\theta_0, \eta_0(\cdot;\theta_0)) -
\dot{\Psi}_2(\theta_0,
\eta_0(\cdot;\theta_0))[\dot{\eta}_0(\cdot;\theta_0)] \Big\}^{-1}
\nonumber \\
&& \quad \ \times \ n^{1/2}\Big\{ (\Psi_n - \Psi)(\theta_0,
\eta_0(\cdot;\theta_0)) \\
&& \qquad \qquad + \ \dot{\Psi}_2(\theta_0,
\eta_0(\cdot;\theta_0))[(\hat{\eta}_n - \eta_0)(\cdot;\theta_0)]
\Big\} + o_{p^*}(1) \, . \hspace{5mm} \mbox{$\Box$}
\end{eqnarray*}

\end{document}